%cannibals.tex: a Plain TeX file George Spahn and Doron Zeilberger
%Variations on the Missionaries and Cannibals Puzzle
%begin macros

\baselineskip=14pt
\parskip=10pt
\def\halmos{\hbox{\vrule height0.15cm width0.01cm\vbox{\hrule height
  0.01cm width0.2cm \vskip0.15cm \hrule height 0.01cm width0.2cm}\vrule
  height0.15cm width 0.01cm}}

\magnification=\magstephalf

\def\1{{\overline{1}}}
\def\2{{\overline{2}}}
\parindent=0pt
\overfullrule=0in

\def\frac#1#2{{#1 \over #2}}
%\headline={\rm  \ifodd\pageno  \RightHead  \else  \LeftHead  \fi}
%\def\RightHead{\centerline{
%Title
%}}
%\def\LeftHead{ \centerline{Doron Zeilberger}}
%end macros
\centerline
{\bf 
Variations on the Missionaries and Cannibals Problem
}
\bigskip
\centerline
{\it George SPAHN and Doron ZEILBERGER}
\centerline
\qquad \qquad \qquad 
{\it Dedicated to Harry Dym (b. Jan. 26, 1938) on  his forthcoming eighty-fifth birthday}
\bigskip

{\bf Abstract}:  We explore both automated and human approaches to the generalized Missionaries and Cannibals problem.

{\bf A Classic Riddle}

Many of us have heard, and quite a few of us solved, the following famous puzzle ([Wi1], and references thereof), that goes back (at least!) to Tartaglia (of cubic equation fame).

{\it
 Three missionaries and three cannibals must cross a river using a boat which can carry at most two people, under the constraint that, 
for both banks, if there are missionaries present on the bank, they cannot be outnumbered by cannibals (if they were, the cannibals would eat the missionaries). 
The boat cannot cross the river by itself with no people on board.
}

This brainteaser was used as a challenge at the very early days of AI (see, e.g. [R], p. 51), and was also a {\it toy model}  in Engineering Design, see the wonderful book [D] (pp. 139-143) by the [biological] brother of the 
academic father of DZ (and hence the academic grandfather of GS).

As indicated in [Wi1] (see also [FCD] and [PS]), in order to solve it {\it nowadays}, one sets up a {\bf directed graph} with the vertices labeled by triples of  non-negative integers $[m,c,b]$ where 

$\bullet$ $m$ is the number of missionaries currently at the first bank

$\bullet$ $c$ is the number of cannibals  currently at the first bank

$\bullet$ $b=1$ if the boat is currently at the first bank, and $b=0$ if it is at the second bank.

It follows that the numbers of missionaries and cannibals in the second bank are $3-m$ and $3-c$ respectively.

In order for no missionary to be eaten, we need at all times:

$\bullet$ If $m>0$ then $m \geq c$

$\bullet$ If $3-m>0$ then $3-m \geq 3-c$, or equivalently if $m<3$ then $c \geq m$.

The edges are

$\bullet$ $[m,c,1] \rightarrow [m-e_1,c-e_2,0] $

$\bullet$ $[m,c,0] \rightarrow [m+e_1,c+e_2,1] $

where $0<e_1+e_2 \leq 2$, and both vertices are legal, where the boat carries $e_1$ missionaries and $e_2$ cannibals.

To solve the puzzle all we need is find a shortest path from the initial state $[3,3,1]$ to the final state $[0,0,0]$.
This can be done with the famous Dijkstra algorithm [Wi2], that we adapted to the present problem.

{\bf  Let the Computer Do it}

Now that we have such high-level computer languages (Maple in our case), it is an easy programming exercise to solve not only the original puzzle, but the following
general version, for any {\bf specific, numerical} {\bf parameters} $M$, $C$, $B$, and $d$, and not only find {\bf one} solution but {\bf all of them}.

So we set up to write Maple code for finding {\bf all} solutions for any inputted $M,C,B,d$, to the following family of puzzles.

{\it
 $M$ missionaries and $C$ cannibals must cross a river using a boat which can carry at most $B$ people, under the constraint that, 
for {\bf both} banks, and for the boat, if there are missionaries and cannibals present there, the number of missionaries must exceed the number of cannibals but at least $d$.
}

[if $d>0$ this means that the cannibals are stronger than the missionaries, so one needs a higher `safety margin'].

Now the number of missionaries and cannibals in the second bank are $M-m$ and $C-c$ respectively.

In order for no missionary to be eaten, we need at all times

$\bullet$ If $m>0$ and $c>0$ then $m - c \geq d$

$\bullet$ If $M-m>0$  and $C-c>0$ then $(M-m) - (C-c) \geq d$

The edges are:

$\bullet$ $[m,c,1] \rightarrow [m-e_1,c-e_2,0] $

$\bullet$ $[m,c,0] \rightarrow [m+e_1,c+e_2,1] $

where $0<e_1+e_2 \leq B$ and both vertices are legal, and $(e_1,e_2)$ satisfies the following condition:

$\bullet$ If $e_1>0$ and $e_2>0$, then $e_1-e_2 \geq d$.

This is implemented in procedure

{\tt Sols(M,C,B,d)} \quad ,

in our Maple package {\tt Cannibals.txt}, obtainable from the front of this article

{\tt https://sites.math.rutgers.edu/\~{}zeilberg/mamarim/mamarimhtml/cannibals.html} \quad.

In particular, to get all four solutions for the original puzzle, type

{\tt Sols(3,3,2,0);} \quad ,

and get in $0.04$ seconds the following four solutions:

$\bullet$ {\tt [[3, 3, 1], [2, 2, 0], [3, 2, 1], [3, 0, 0], [3, 1, 1], [1, 1, 0], [2, 2, 1],[0, 2, 0], [0, 3, 1], [0, 1, 0], [0, 2, 1], [0, 0, 0]]} \quad ;

$\bullet$ {\tt [[3, 3, 1], [2, 2, 0],[3, 2, 1], [3, 0, 0], [3, 1, 1], [1, 1, 0], [2, 2, 1], [0, 2, 0], [0, 3, 1], [0, 1, 0], [1, 1, 1], [0, 0, 0]]} \quad ;

$\bullet$ {\tt [[3, 3, 1], [3, 1, 0], [3, 2, 1], [3, 0, 0], [3, 1, 1], [1, 1, 0], [2, 2, 1], [0, 2, 0], [0, 3, 1], [0, 1, 0], [0, 2, 1], [0,0, 0]]} \quad ;

$\bullet$ {\tt [[3, 3, 1], [3, 1, 0], [3, 2, 1], [3, 0, 0], [3, 1, 1], [1, 1, 0], [2,2, 1], [0, 2, 0], [0, 3, 1], [0, 1, 0], [1, 1, 1], [0, 0, 0]]} \quad .

For any solution, {\tt S}, procedure {\tt SO(M,C,B,d,S)}, in our Maple package {\tt Cannibals.txt},  spells  out the solution, very verbosely. See the output file

{\tt https://sites.math.rutgers.edu/\~{}zeilberg/tokhniot/oCannibals1.txt} \quad .

To see all $25$ solutions with $5$ missionaries and $5$ cannibals, boat-capacity $3$ and safety-margin still $0$, as well as the spelled-out 
solution for one of the solutions see

{\tt https://sites.math.rutgers.edu/\~{}zeilberg/tokhniot/oCannibals2.txt} \quad .

The above-mentioned front of our article contains a few more sample solutions, but readers (who have access to Maple) are 
welcome to find solutions to  many other cases.

{\bf  Paths in General Directed Graphs and Linear Algebra}

After we completed the above Maple package, {\tt Cannibals.txt}, that focuses on one specific class of problems, we realized that we may
just as well treat the general case of finding all walks in an {\it arbitrary} directed graph, from one designated vertex to another.
So we wrote a much more general  Maple package, {\tt DiGpaths.txt}, available from

{\tt https://sites.math.rutgers.edu/\~{}zeilberg/tokhniot/DiGpaths.txt} \quad .

In order to represent a directed graph in the computer, we used the data structure {\bf list of sets}
$$
[N_1, ..., N_n] \quad .
$$
Here we assumed that the set of vertices is  $\{1, \dots , n\}$, and  for each vertex $1 \leq i \leq n$,
$N_i$ is the set of vertices $j$ such that there is a directed edge from vertex $i$  to vertex $j$. 
In other words the {\it set of outgoing neighbors of $i$}.

We also make the convention that the originating vertex is $1$ and the terminal vertex is $n$. Our function {\tt RG(n,p)} can construct
many random examples. We want to find the set of shortest paths from $1$ to $n$. In order to do that, we
recursively construct longer and longer {\it self-avoiding} paths until we encounter the terminal vertex $n$.
If after $n$ steps we did not reach it, we declare failure.

For any specific class of puzzles, e.g. Missionaries and Cannibals with general parameters, 
or  other `river-crossing' riddles,
we construct the relevant directed graph. See procedure {\tt MCgraph(M,C,B,d)} in {\tt DiGpaths.txt}.

If we only want to find the {\bf number} of such paths then we can use linear algebra. Recall that
the {\bf adjacency} matrix of a directed graph $G$ is the $0-1$ matrix defined by
$$
A_{ij} \, = \, \cases{1, & if there is a directed edge from $i$ to $j$; \cr
                      0, &otherwise.\cr}
$$

The number of paths of length $k$ between vertex $i$ and vertex $j$ is the $(i,j)$ entry of $A^k$. 
Since we are only interested in {\it shortest} paths, their length is at most $n$, hence in order
to find the length of the shortest paths, and their number, we keep raising the power $A^k$
until, for the first time, the $(ij)$ entry is non-zero. That first successful $k$ would give us the 
desired {\bf minimal length}, and the corresponding $(ij)$ entry is the desired {\bf number} of  shortest paths
(i.e. the number of solutions to our puzzle).
Note that this is {\it number-crunching}, and it only tells you whether a solution {\bf exists}, and if it does, 
{\bf the number of solutions}, but it does {\bf not} tell you how to actually solve the puzzle, i.e. to find
at least one legal path in the underlying graph.

This is implemented in procedure {\tt NuPaths} of the Maple package {\tt DiGpaths.txt} .

Since Maple can do {\bf symbol-crunching} so well, we can use {\bf linear algebra} not just to find the {\bf number} of paths,
but to actually find {\bf all of them}. Define the {\bf symbolic adjacency matrix} by
$$
S_{ij} \, = \, \cases{a_{ij}, & if there is a directed edge from $i$ to $j$; \cr
                      0, &otherwise.\cr}
$$
Then, since we already know, from the numeric adjacency matrix, what is the smallest $k$ such that
the $(i,j)$ entry of $A^k$ is non-zero, and that means that
for that very same $k$,   $S^k$ is also non-zero, that $(i,j)$ entry of $S^k$ gives a certain {\bf polynomial} in
the indeterminates $a_{ij}$  (where $ij$ is an edge) where each monomial corresponds to a path that is easy to
reconstruct, since we are given the set of participating edges in each such path.

This is implemented in procedure {\tt WtPath} in the Maple package {\tt DiGpaths.txt}. Alas, to our disappointment, for larger graphs, this
took much longer than the straightforward approach of generating all paths, mentioned above (procedure {\tt PathsG(G)} in the Maple package {\tt DiGpaths.txt}, and {\tt SolveMC(M,C,B,d)} when applied
to  Missionaries and Cannibals problems).

One can also generate interesting sequences for `infinite families' of Missionaries and Cannibals puzzles. Procedure
{\tt SeqrBd(r,B,d,K)} outputs the first $K$ terms of the sequence enumerating the
number of solutions to the Missionaries and Cannibals problem with 

$i$ missionaries, $i+r$ cannibals, boat side $B$ and safety margin $d$. 

See the output file

{\tt https://sites.math.rutgers.edu/\~{}zeilberg/tokhniot/oDiGpaths2.txt} \quad.

In particular  {\tt SeqrBd(5,3,1,8)} outputs

$$
[4, 4, 13, 21, 34, 55,89, 144] \quad .
$$

This leads to the following conjecture.

{\bf Conjecture}: The number of solutions  to the Missionaries and Cannibals problem
where there are  $i + 5$ missionaries, $i$  cannibals , the Boat capacity is $3$  and the safety margin is $1$ is the Fibonacci number $F_{i+4}$, for $i \geq 3$.

In the next section we will describe an even more efficient approach for generating such conjectures.

Sometimes the enumerating sequences are eventually constant. The next proposition was first conjectured experimentally, using our Maple packages, and later we found
a {\it human-generated proof}.

{\bf Proposition}: The number of solutions with $i$ missionaries and $i$ cannibals, boat size $4$, and safety margin $0$,
is always $361$, for all $i \geq 7$.

{\bf Human-Generated Proof:}
If $n$ is sufficiently large compared to $B$, the boat size, then at some point during the solution there will need to be many missionaries and cannibals on both sides of the river. When this happens, the only legal moves are sending across boats that contain an equal number of missionaries and cannibals. If the missionaries exceed the cannibals in a location, that implies that the cannibals exceed the missionaries in another location, which is an illegal state.

Let $(a,b)$ denote a move that sends $a$ missionaries and $b$ cannibals across the river. Let $-(a,b)$ be a similar moves that brings people in the other direction. If the boat size is less than 4, and there are more than 5 missionaries and cannibals on both sides of the river, the only legal moves will be $\pm(1,1)$. It is clear that no progress can be made with only these moves. This shows that for $B < 4$, and sufficiently large $n$, there will be no solutions. It is not too hard to show that for $B=2$ there exist solutions only when $n < 4$, and that for $B=3$ there exist solutions only when $n < 6$. 

When $B=4$, we now can make progress for large $n$. We can alternate moves $(2,2)$ and $-(1,1)$ to slowly ferry the people across without violating the rules. Using this approach we can solve the puzzle in $2n-3$ moves. It is clear that during the middle of a solution, it is not possible to be more efficient. But maybe we can be more efficient at the start or end, by exploiting the fact that the cannibals are allowed outnumber the missionaries if either are not present on a side of the river. Let's say the start of the solution is the moves before 4 missionaries and 4 cannibals have crossed. Once 4 and 4 have crossed, we are now restricted to alternating  $(2,2)$ and $-(1,1)$. After the first two moves, at most 3 people will have crossed. Thus it is clear that it is not possible to do the start in 3 moves. Thus the minimum possible number of moves for the start of the solution is 5 moves. It turns out there are exactly 19 ways to perform the start of the puzzle in 5 moves. Symmetrically, there are also 19 ways to perform the end of a solution as efficiently as possible. For $n \geq 7$, the start of the solution and the end of solution do not intersect, so we conclude that for $n \geq 7$, the number of solutions that use the fewest numbers of moves will be $19 \cdot 19 =361$.  \halmos

{\bf  Using ``High School Algebra" for Efficient Enumeration}

While finding {\it one} solution (or path) can be done fast, using   Dijkstra's algorithm, finding {\it all of them} takes much longer, since there are exponentially many of them.
Of course we can do it as in the previous section, transcribing the puzzle into a directed graph, finding the adjacency matrix, etc. but the resulting
graph has special structure, and it would be nice to take advantage of this.

We will now describe a naive `high-school algebra' approach for quickly counting the number of solutions (as opposed to finding all of them, or even one of them).
In particular, if we find out, for a particular choice of parameters (for example, $(M,C,B,d)=(4,4,2,0)$), that the number of solutions is $0$, we would know
that there is {\bf no way} of solving the puzzle.

Let's consider a more general river-crossing puzzle, that includes the even-more-famous-and-older Cabbage-Sheep-Wolf puzzle [Wi3].

There are $k$ species, let's call them $1, \dots, k$, and there are $A_1,A_2, \dots , A_k$ individuals from each kind. 
If at any time, there are $a_1, \dots, a_k$ individuals of species $1, \dots, k$ respectively,
at the starting bank, then the vectors $[a_1,a_2, \dots, a_k]$
must satisfy a set of conditions $C(a_1, \dots, a_k)$. There is another set of conditions $C'(b_1, \dots, b_k)$ regarding
what the boat can carry.

Let's introduce $k$ `formal' variables (alias {\it indeterminates}) $x_1, \dots. x_k$ and define the {\bf crossing polynomial}
$$
P(x_1, \dots, x_k) := \sum_{(b_1, \dots, b_k)} x_1^{b_1} \cdots x_k^{b_k} \quad,
$$
where the sum is over all vectors of non-positive integers $(b_1, \dots, b_k)$ satisfying the conditions $C'(b_1, \dots, b_k)$.

For example, for the original puzzle where $x_1$ is the variable corresponding to missionaries and $x_2$ is the variable corresponding to cannibals 
(and the boat size is $2$)
$$
P(x_1,x_2)= x_1+x_1^2 +x_2+x_2^2 + x_1 x_2 \quad .
$$

Let's also introduce a ``clean-up" linear operator $T$ that is defined on {\bf poly}nomials  in the variables $x_1, \dots, x_k$, by first defining it on
{\bf mono}mials $x_1^{a_1} \cdots x_k^{a_k}$ by
$$
T(x_1^{a_1} \cdots x_k^{a_k}) \, = \, \cases{ x_1^{a_1} \cdots x_k^{a_k}  , & if $C(a_1, \dots, a_k)$ is true ; \cr
                      0, &otherwise,\cr} 
$$
and extending it {\bf linearly}.

We start with the monomial representing the {\bf initial position} $x_1^{A_1} \cdots x_k^{A_k}$. Going in the forward direction corresponds
to {\bf multiplying} by
$$
P(x_1^{-1}, \dots, x_k^{-1}) \quad,
$$
since you reduce the population of the first bank, while going  back to the starting bank corresponds to multiplying by
$$
P(x_1, \dots, x_k) \quad ,
$$
since you increase the population of the first bank.

At any stage, after multiplying by $P(x_1, \dots, x_k)$ or $P(x_1^{-1}, \dots, x_k^{-1})$, we have to ``clean up'' by applying the operator $T$.

This leads to the following sequence of polynomials
$$
f_0(x_1, \dots, x_k) = x_1^{A_1} \cdots x_k^{A_k} \quad,
$$
and for $i=1,2,\dots$
$$
g_i(x_1, \dots, x_k) = T[\,P(x_1^{-1}, \dots, x_k^{-1}) f_{i-1}(x_1, \cdots, x_k)\,]
$$
$$
f_i(x_1, \dots, x_k) = T[\,P(x_1, \dots, x_k) g_{i}(x_1, \cdots, x_k)\,] \quad .
$$

We keep going until we reach an $i$ where $g_i$ has {\bf non-zero} constant term. If we do reach such an $i$ we know that the
puzzle can be solved in $2i-1$ crossing ($(i-1)$ double crossings followed by the final one), and the constant term of that
lucky $g_i$ is the number of solutions to the puzzle (equivalently the number of shortest paths). Of course
it is very possible that we will never reach that state.

Since any shortest path must be {\it self-avoiding}, we have an upper bound for the length of a shortest path, if it exists.
Let's call the number of states $T$. If by the time $i=T+1$ none of the $g_i$ have a non-zero constant term, we find out
that there are no solutions.

Let's illustrate our method with the original puzzle of $3$ missionaries, $3$ cannibals, boat size $2$, and safety margin $0$.
Recall that for this case $P(x_1,x_2)=x_1+x_1^2 +x_2+x_2^2 + x 1 x_2$.
$$
f_0=x_1^3\, x_2^3 \quad.
$$
$$
f_0 \cdot P(x_1^{-1},x_2^{-1})= 
(x_1^3\, x_2^3) \cdot (x_1^{-1}+x_1^{-2} +x_2^{-1}+x_2^{-2} + x_1^{-1} x_2^{-1})
=x_{1}^{3} x_{2}^{2}+x_{1}^{2} x_{2}^{3}+x_{1}^{3} x_{2}+x_{2}^{2} x_{1}^{2}+x_{1} x_{2}^{3}
\quad.
$$
``Cleaning up", i.e. applying $T$, we get
$$
g_1(x_1,x_2)=x_{1}^{3} x_{2}^{2}+x_{1}^{3} x_{2}+x_{2}^{2} x_{1}^{2} \quad.
$$
Now going back to the other bank
$$
g_1(x_1,x_2) \cdot  P(x_1,x_2) \, = \,
(x_{1}^{3} x_{2}^{2}+x_{1}^{3} x_{2}+x_{2}^{2} x_{1}^{2})( x_1+x_1^2 +x_2+x_2^{2} + x_1 x_2)=
$$
$$
x_{1}^{5} x_{2}^{2}+x_{1}^{4} x_{2}^{3}+x_{1}^{3} x_{2}^{4}+x_{1}^{5} x_{2}+3 x_{1}^{4} x_{2}^{2}+3 x_{1}^{3} x_{2}^{3}+x_{1}^{2} x_{2}^{4}+x_{1}^{4} x_{2}+2 x_{1}^{3} x_{2}^{2}+x_{1}^{2} x_{2}^{3}
 \quad.
$$
``Cleaning up", i.e. applying $T$, by discarding all `illegal monomials', we get
$$
f_1(x_1,x_2)= 3 x_{1}^{3} x_{2}^{3}+2 x_{1}^{3} x_{2}^{2} \quad.
$$
Moving right along
$$
g_2(x_1,x_2)= 3 x_{1}^{3} x_{2}^{2}+5 x_{1}^{3} x_{2}+5 x_{2}^{2} x_{1}^{2}+2 x_{1}^{3} \quad ,
$$
$$
f_2(x_1,x_2)= 13 x_{1}^{3} x_{2}^{3}+12 x_{1}^{3} x_{2}^{2}+2 x_{1}^{3} x_{2} \quad,
$$
$$
g_3(x_1,x_2)=  13 x_{1}^{3} x_{2}^{2}+25 x_{1}^{3} x_{2}+25 x_{2}^{2} x_{1}^{2}+14 x_{1}^{3}+2 x_{1} x_{2} \quad,
$$
$$
f_3(x_1,x_2)= 63 x_{1}^{3} x_{2}^{3}+64 x_{1}^{3} x_{2}^{2}+16 x_{1}^{3} x_{2}+2 x_{2}^{2} x_{1}^{2} \quad,
$$
$$
g_4(x_1,x_2)=63 x_{1}^{3} x_{2}^{2}+127 x_{1}^{3} x_{2}+127 x_{2}^{2} x_{1}^{2}+80 x_{1}^{3}+18 x_{1} x_{2}+2 x_{2}^{2} \quad,
$$
$$
f_4(x_1,x_2)=317 x_{1}^{3} x_{2}^{3}+334 x_{1}^{3} x_{2}^{2}+98 x_{1}^{3} x_{2}+20 x_{2}^{2} x_{1}^{2}+2 x_{2}^{3} \quad,
$$
$$
g_5(x_1,x_2)= 317 x_{1}^{3} x_{2}^{2}+651 x_{1}^{3} x_{2}+651 x_{2}^{2} x_{1}^{2}+432 x_{1}^{3}+118 x_{1} x_{2}+22 x_{2}^{2}+2 x_{2} \quad,
$$
$$
f_5(x_1,x_2)= 1619 x_{1}^{3} x_{2}^{3}+1734 x_{1}^{3} x_{2}^{2}+550 x_{1}^{3} x_{2}+140 x_{2}^{2} x_{1}^{2}+24 x_{2}^{3}+2 x_{1} x_{2}+2 x_{2}^{2} \quad,
$$
and {\bf finally}:
$$
g_6(x_1,x_2)= 619 x_{1}^{3} x_{2}^{2}+3353 x_{1}^{3} x_{2}+3353 x_{2}^{2} x_{1}^{2}+2284 x_{1}^{3}+690 x_{1} x_{2}+164 x_{2}^{2}+28 x_{2}+4 \quad .
$$
Success! We have found that with $i=6$ $g_i(x_1,x_2)$ has a non-zero constant term, that happens to be $4$.
Hence the puzzle can be solved with $2 \cdot 6-1=11$ crossings, and the number of solutions is $4$.

If we try to solve the puzzle with $4$ missionaries and $4$ cannibals, we get that the $g_i$ have no non-zero constant term for all $i \leq 14$, and since there are $13$ legal states,
it means that there are no solutions.

Since Maple is so good with high-school algebra, this is very fast, and we can investigate many infinite families,  getting similar
conjectures to the Fibonacci conjecture from the previous section. To enjoy $139$ such propositions, see

{\tt https://sites.math.rutgers.edu/\~{}zeilberg/tokhniot/oRiverCrossing2.txt} \quad.

Let's just reproduce one of these conjectures, that we sure are theorems.

{\bf Proposition}: Let $a(n)$ be the number of ways of safely transporting $n+9$ missionaries, $n$ cannibals, with a boat that can have at most
two passengers, and such that at no time, at either bank,  should the number of cannibals exceed the number of missionaries (if there are missionaries present), then
$$
\sum_{n=0}^{\infty} a(n) \, x^n \,= \,
$$
$$
-\frac{1774224 x^{7}-63279616 x^{6}-54735368 x^{5}+31754164 x^{4}-2667061 x^{3}-736742 x^{2}-6726 x -1}{4 x^{4}-384 x^{3}+337 x^{2}-39 x +1} \quad .
$$

{\bf Let Humans do some General Thinking}

We now consider the problem of determining for which values of the parameters there exists a solution. Again we have:

   $\bullet$ M, the number of missionaries, should be at least 1
   
    $\bullet$ C, the number of cannibals, should be at least 1
    
   $\bullet$ B, the boat size, should be at least 2 
   
    $\bullet$ d, if both cannibals and missionaries are present (on either side of the river or in the boat) then the number of missionaries must be at least the number of cannibals plus d.
 \halmos

Tackling the problem {\it in general} is hard, but there are a lot of cases where it is easy to show that a solution exists. We do so by demonstrating a strategy that illustrates how to solve the puzzle in a specific case, and then mention the requirements on the parameters for that strategy to be successful. The strategies provide no insight into what happens when their requirements are not met, so these strategies give a set of sufficient conditions for when a solutions exists, but do not provide necessary conditions.

{\bf ``Two Boat Strategy":}
This strategy only requires two people to be in the boat at any given time (hence the name), and therefore works for any value of $B \geq 2$. To start, we can send a single missionary across by sending 2 missionaries across in the boat, and then having one come back with the boat. Call this sequence of 2 moves P. By repeating P, we can send missionaries across until P is no longer legal. We can send over $M - C - d - 1$ missionaries in this way. Now consider the same operation, but with cannibals instead of missionaries. Call this operation Q. If Q is legal, we can do Q, and then complete the puzzle by alternating P and Q until all the cannibals are across.  The requirement for Q to be legal is that \quad $(M - C - d - 1) \geq 2 + d $ \quad which simplifies to:

$$ M-C \geq 2d+3$$

We now know that if the above condition is satisfied, then the puzzle is solvable.

{\bf ``Big Boat Strategy \#1":} If the boat is big enough relative to the number of Cannibals, we can send all of the missionaries across before sending any cannibals across. Similar to Two Boat, we start by sending $M - C - d - 1$ missionaries across and bringing the boat back. The number of missionaries remaining is then $C + d + 1$. If $$ B \geq C + d + 1 $$ we can proceed to send the rest of the missionaries over, send $B-1$ back, use them to ferry a single cannibal across, send the cannibal back, and then the cannibals can all ferry themselves across.

{\bf ``Big Boat Strategy \#2":}
This case is very trivial, but does require its own approach. If the boat is so big that it can carry all the missionaries in it, then we can use this strategy. Send 2 cannibals across, send 1 back, send all the missionaries over, send 1 cannibal back, send the rest of the cannibals across.
$$ B \geq M, \quad C \geq 2$$

{\bf ``Split Cannibals Strategy":}
In this strategy we send half of the cannibals over, and then all the missionaries over. It is similar to Big Boat in that it only works if the Boat is large compared to half the number of cannibals. The strategy is slightly different depending on if there is an even number or odd number of cannibals. 

First consider the even case. We send over half the cannibals. Then we can send a boat of $C/2 + d + 1$ missionaries. This requires $$ B > C/2 + d + 1$$ and $$ M - (C/2 + d + 1) \geq C/2 + d$$ which simplifies to $$ M - C \geq 2d+1$$ Then the missionaries can ferry themselves across, and then the cannibals can ferry themselves across.

In the odd case, we send $\lceil C/2 \rceil$ cannibals across, then $\lceil C/2 \rceil + d + 1$ missionaries across, and we can proceed similarly if $$ B > \lceil C/2 \rceil + d + 1$$ and $$ M - (\lceil C/2 \rceil + d + 1) \geq \lfloor C/2 \rfloor + d$$ which simplifies to $$ M - C \geq 2d+1$$
However, it is still sometimes possible if $$ B = \lfloor C/2 \rfloor + d + 1$$ 
This is an open problem.

We can combine the even and odd cases and conclude that the strategy can be applied when 
$$ 
M - C \geq 2d+1 \quad AND \quad  B > \lceil C/2 \rceil + d + 1
$$

{\bf ``Simultaneous Ferry Strategy"}:

Here we start by sending $d$ missionaries over, and then we repeatedly send $d+1$ missionaries and $1$ cannibal over, and $d$ missionaries back. Since cannibals and missionaries are in all 3 places at once, this requires 
$$ M - C \geq 3d$$ in addition to $$ B \geq d+2$$

We now summarize the conditions required for each strategy to be applied.

$$
\vbox
{\tabskip=0pt
\halign{\tabskip=4cm #\hfil&#\hfil\tabskip=0pt \cr
Two Boat &  $ M-C \geq 2d+3$\cr
Big Boat 1 &$ B \geq C + d + 1 $\cr
Big Boat 2 &$ B \geq M$ \quad AND \quad $ C \geq 2$\cr
Split Can &$ M - C \geq 2d+1$ \quad AND \quad $ B > \lceil C/2 \rceil + d + 1$ \cr
Simul &$ M - C \geq 3d$ \quad AND \quad $ B \geq d+2$\cr
}
}
$$

Also note that the condition that $M-C\geq d$ is assumed.

The above strategies do not say too much about what happens when $d=0$. Since this was the original version of the puzzle it deserves some extra attention. It turns out that we can do much better than the listed strategies, and in almost all cases it is very simple to do so. We consider if there is some slack, that is $M > C$.

{\bf ``d=0,  M $>$ C": } Send over a boat with 1 missionary and 1 cannibal. Send 1 cannibal back. Call this operation Q. Repeat Q until all the cannibals are across, then have the missionaries ferry the remaining missionaries across. Done.

When $M=C$ and $B \geq 4$, another simple strategy applies.

{\bf ``d=0,  M =C, B $\geq$ 4":} Send over 2 missionaries and 2 cannibals. Send 1 missionary and 1 cannibal back. Repeat until everyone is across. 

The remaining cases are when $M=C$, and either $B=2$ or $B=3$. For $B=2$ it turns out to be doable if and only if $M = C \leq 3$ and for $B=3$ it turns out to be doable if and only if $M = C \leq 5$. These cases are not easily described by a simple strategy and make for a fun puzzle.

{\bf Necessary conditions?}
Are the above conditions necessary for the riddle to be solved? Not entirely. One example is that sometimes the Split cannibals strategy can still be applied when the condition is not satisfied. With more computation power, we should get a better idea of any potential cases where a strategy not described can be used.

{\bf References}

[D] Clive L. Dym, {\it ``Engineering Design, A Synthesis of Views''}, Cambridge University Press, 1994. 

[FCD]  Robert Fraley,  Kenneth L. Cooke, and Peter  Detrick, {\it Graphical Solution of Difficult Crossing Puzzles}. Mathematics Magazine {\bf 39}(3) (May 1966).  151-157.

[PS] Ian Pressman and David Singmaster, {\it The jealous husbands and The missionaries and cannibals}, Mathematical Gazette {\bf 73}, No. 464, (June 1989), 73-81.

[R] Elaine Rich,{\it ``Artificial Intelligence''}, McGraw Hill, 1983.

[Wi1] Wikipedia, {\it Missionaries and cannibals problem}.

[Wi2] Wikipedia, {\it  Dijkstra's algorithm}.

[Wi3] Wikipedia, {\it  Wolf, goat and cabbage problem}

\bigskip
\hrule
\bigskip
George Spahn and Doron Zeilberger, Department of Mathematics, Rutgers University (New Brunswick), Hill Center-Busch Campus, 110 Frelinghuysen
Rd., Piscataway, NJ 08854-8019, USA. \hfill\break
Email: {\tt  gs828 at math dot rutgers dot edu} \quad, \quad {\tt DoronZeil] at gmail dot com}   \quad .

Written: {\bf Oct. 21, 2022}.

\end